\documentstyle[11pt, leqno, epsfig]{article}
\input amssym.def
\title{Eigenvalue estimates of the Dirac operator depending on the Ricci 
tensor.  
\footnote{Supported by the SFB 288 of the DFG.}}
\author{Thomas Friedrich and Klaus-Dieter Kirchberg (Berlin)}
\date{\today}

\begin{document}

\newcommand{\D}{\displaystyle}
\newcommand{\upsp}{\phantom{l}}
\newcommand{\downsp}{\phantom{q}}

\maketitle

\mbox{} \hrulefill \mbox{}\\

\newcommand{\vol}{\mbox{vol} \, }
\newcommand{\grad}{\mbox{grad} \, }
\newcommand{\sneun}{Spin(9)}

\begin{abstract} We prove a new lower bound for the first eigenvalue of the 
Dirac operator on a compact Riemannian spin manifold by refined Weitzenb\"ock
techniques. It applies to manifolds with harmonic curvature tensor and depends on the Ricci tensor. Examples show how it behaves compared to other known 
bounds.
\end{abstract}

{\small
{\it Subj. Class.:} Differential Geometry.
}\\

\mbox{} \hrulefill \mbox{}\\

\setcounter{section}{-1}


\section{Introduction}

If $M^n$ is a compact Riemannian spin manifold with 
positive scalar curvature $R$, then each eigenvalue $\lambda$ of the Dirac 
operator $D$ satisfies the inequality
\begin{equation}
\lambda^2 \ge \frac{n R_0}{4(n-1)} \ ,
\end{equation}

where $R_0$ is the minimum of $R$ on $M^n$. The estimate $(1)$ is sharp
in the sense that there exist manifolds for which $(1)$ is an equality for the 
first eigenvalue $\lambda_1$ of $D$. If this is the case, then each
eigenspinor $\psi$ corresponding to $\lambda_1$ is a Killing spinor with the
Killing number $\lambda_1/n$, i.e., $\psi$ is a solution of the field equation
\begin{equation}
\nabla_X \psi + \frac{\lambda_1}{n} X \cdot \psi =0
\end{equation}

and $M^n$ must be an Einstein space (see \cite{F1}). A generalization of 
this inequality was proved in the paper 
\cite{Hi}, where a conformal lower bound for the spectrum of the Dirac 
operator occured. Moreover, for special Riemannian manifolds better estimates 
for the eigenvalues of the Dirac operator are known, see \cite{Ki}, \cite{KSW}. However, all these estimates of the spectrum of the Dirac operators depend 
only on the scalar curvature of the underlying manifold. Therefore it is a 
natural question whether or not one may relate the spectrum of the Dirac 
operator to more refined curvature data.\\

In this paper we shall prove an estimate depending on the Ricci tensor for the 
eigenvalues of the Dirac operator on compact Riemannian manifolds with harmonic
curvature tensor. The main idea is the investigation of the differential 
operators 
\begin{displaymath}
Q^t : \Gamma (S) \to \Gamma (TM^n \otimes S) 
\end{displaymath}
depending on a real parameter $t \in {\Bbb R}$ and defined by
\begin{displaymath}{\textstyle
Q^t_X \psi := \nabla_X D \psi + \frac{1}{n} X \cdot D^2 \psi + t \cdot
(\mbox{Ric} - \frac{R}{n})(X) \cdot \psi ,} 
\end{displaymath}

where $\mbox{Ric}$ denotes the Ricci tensor. Under the assumption that the 
curvature tensor is harmonic we prove a formula
expressing the length $|Q^t\psi|^2$ by the Dirac operator $D\psi$, the
covariant derivatives $\nabla D\psi$ and $\nabla \psi$ as well as by some
curvature terms (Theorem 1.6). Integrating this formula we obtain, for any 
$t \ge 0$, an inequality for the eigenvalues of the Dirac operator
depending on the scalar curvature, the minimum of the eigenvalues of the
Ricci tensor and its length. An optimal choice of the parameter $t$ bounds
the spectrum of the Dirac operator from below. For example, we prove the 
inequality
\begin{displaymath}
\lambda^2 > \frac{1}{4} \cdot \frac{| \mbox{Ric}|^2_0}{| \mbox{Ric}|_0 
\sqrt{\frac{n-1}{n}} + |\kappa_0|} 
\end{displaymath}
for compact Riemannian spin manifolds with harmonic curvature tensor and 
vanishing scalar curvature, where $\kappa_0$ and $|\mbox{Ric}|_0$ denote the minimum of the eigenvalues and the length of the Ricci tensor, respectively.

\section{The Weitzenb\"ock formula for the operator $Q^t$}

First of all let us fix some notations. In the following $(X_1, \ldots,
X_n)$ is always any local frame of vector fields, and $(X^1, \ldots , 
X^n)$ is the associated frame defined by $X^k := g^{kl}X_l$, where the  
$g^{kl}$ denote the components of the inverse of the Riemannian metric 
$(g_{kl}):=(g(X_k, X_l))$. Using the twistor operator 
(see \cite{BFGK}, Section 1.4)
\begin{displaymath}
{\cal D}: \Gamma (S) \longrightarrow \Gamma (TM^n \otimes S)
\end{displaymath}
locally given by ${\cal D}{\psi} := X^k \otimes {\cal D}_{X_k} \psi$ and
${\cal D}_X \psi \ := \ \nabla_X\psi + \frac{1}{n} X \cdot D \psi$, we may
rewrite the operator $Q^t$ as

\begin{equation}{\textstyle
Q^t \psi \ = \ {\cal D} D  \psi + t \cdot X^k \otimes ( \mbox{Ric} - 
\frac{R}{n} )(X_k) \cdot \psi .}
\end{equation}

The image of the twistor operator ${\cal D}$ is contained in the kernel of the
 Clifford multiplication $\mu :TM^n \otimes S \to S$, i.e.,

\begin{equation}
\mu ({\cal D} \psi)\ = \ X^k \cdot {\cal D}_{X_k} \psi =0 \ .
\end{equation}

As endomorphisms acting on the spinor bundle the following identities are well
known:
\begin{equation}
X^k \cdot \mbox{Ric} (X_k) \ = \  \mbox{Ric} (X_k) \cdot X^k \ = \  - R, \quad
{\textstyle
X^k \cdot ( \mbox{Ric} - \frac{R}{n} ) (X_k) \ = \ 0 \ .} 
\end{equation}

In particular, we see that the image of the operators $Q^t$ is contained in 
the kernel of the Clifford multiplication. By definition, a spinor 
field $\psi$ belongs to the kernel of the operator $Q^t$ if and only if
it satisfies the equation
\begin{equation}{\textstyle
\nabla_X  D \psi + \frac{1}{n} X \cdot D^2 \psi + t \cdot ( \mbox{Ric} - 
\frac{R}{n} ) (X) \cdot \psi \ = \ 0}
\end{equation}

for each vector field $X$.
In the following we shall use the Weitzenb\"ock formula
\begin{equation}{\textstyle
|{\cal D} \psi |^2 \ = \ |\nabla \psi |^2 - \frac{1}{n} |D \psi |^2} 
\end{equation}

for the twistor operator $\cal{D}$.\\

{\bf Lemma 1.1:} {\it For any spinor field $\psi \in \Gamma (S)$, the following formula holds:}
\begin{equation}
\begin{array}{l}
|Q^t \psi |^2 \ = \ | \nabla D \psi |^2 - \frac{1}{n} |D^2 \psi |^2 +2t 
\cdot \frac{R}{n} \cdot \mbox{Re} ( \langle D^2 \psi , \psi \rangle 
)+ t^2 \cdot |\mbox{Ric} - \frac{R}{n} |^2 \cdot |\psi |^2  \\ \mbox{}\\
\mbox{} \quad \quad \quad \quad - 2t \cdot \mbox{Re} \Big(\langle 
\mbox{Ric} (X^k) \cdot \nabla_{X_k} D \psi, \, \psi \rangle \Big)
\end{array} 
\end{equation}

{\bf Proof:} Using the formulas $(3), (4)$ and $(7)$ we have
\begin{displaymath}
\begin{array}{lll} 
|Q^t \psi |^2 &=& \langle Q^t_{X_k} \psi , Q^t_{X^k} \psi \rangle \\[1em]
&\stackrel{(3)}{=}&\Big\langle {\cal D}_{X_k} D \psi + t (\mbox{Ric} - 
\frac{R}{n}) (X_k) \cdot \psi
, {\cal D}_{X^k} D \psi + t ( \mbox{Ric} - \frac{R}{n} ) (X^k) 
\cdot \psi \Big\rangle\\[1em]
&=& |{\cal D}D \psi |^2 - 2t \cdot \mbox{Re} \left( \Big\langle ( 
\mbox{Ric} - \frac{R}{n} ) (X^k)
 {\cal D}_{X_k} D \psi , \psi \Big\rangle \right)+ t^2 \cdot | \mbox{Ric} 
- \frac{R}{n} |^2 \cdot |\psi|^2 \\[1em]
&\stackrel{(4)}{=}& |{\cal D}D \psi |^2 + t^2 \cdot |\mbox{Ric} - 
\frac{R}{n} |^2 \cdot |\psi|^2 - 2t \cdot \mbox{Re} \left( \langle 
\mbox{Ric} (X^k) \cdot {\cal D}_{X_k} D \psi , \psi \rangle \right) \\[1em] 
&\stackrel{(7)}{=}& | \nabla D \psi |^2 - \frac{1}{n} |D^2 \psi |^2 + t^2 
\cdot \Big| \mbox{Ric} - \frac{R}{n} \Big|^2 \cdot |\psi |^2  +2t \cdot 
\frac{R}{n} \cdot \mbox{Re} \left( \langle D^2 \psi , \psi \rangle \right) 
\\[1em]
&&- 2t \cdot \mbox{Re} \left(\langle \mbox{Ric} (X^k) \cdot \nabla_{X_k} D \psi , \psi \rangle \right) . \hfill \Box
\end{array} 
\end{displaymath}

Equation $(8)$ is a preliminary version of the Weitzenb\"ock formula, which we
will apply in the proof of our main result. Our next aim is to express the 
uncontrollable last term  on the right-hand side by terms that are 
controllable. For this purpose we need a condition on the covariant derivative of the 
Ricci tensor. For vector fields $X,Y$, we use the notation
\begin{displaymath}
\nabla_{X,Y} := \nabla_X \nabla_Y - \nabla_{\nabla_X Y}
\end{displaymath}

for the corresponding tensorial derivatives of second order in $TM^n$ as well 
as in $S$. By $K$ we denote the Riemannian curvature tensor and by $C$ the 
curvature tensor in the spinor bundle $S$. Then, for all $X,Y, Z \in \Gamma 
(TM^n)$ and all $\psi \in \Gamma (S)$, we have
\begin{displaymath}
K(X,Y)(Z) \ = \ \nabla_{X,Y} Z - \nabla_{Y,X} Z \quad , 
\quad C(X,Y) \psi \ = \ \nabla_{X,Y} \psi - \nabla_{Y,X} \psi 
\end{displaymath}
as well as the well known relation between the two curvatures
\begin{equation}{\textstyle
C(X,Y) \cdot \psi \, =\, \frac{1}{4}\, X^k \cdot K(X,Y) (X_k) \cdot \psi 
\, = \,  \frac{1}{4} \, g(K(X,Y)(X^k) , X^l) X_k \cdot X_l \cdot \psi \ .}
\end{equation}

Considering $C$ as a map from $\Gamma (S)$ 
to $\Gamma (TM^n \otimes TM^n \otimes S)$ locally defined by
\begin{displaymath}
C \psi \ := \ X^k \otimes X^l \otimes C(X_k, X_l) \psi , 
\end{displaymath}

the length $|C \psi |^2$ is just the scalar product
\begin{displaymath}
|C \psi |^2 = \langle C(X_k, X_l) \psi , C(X^k, X^l ) \psi \rangle . 
\end{displaymath}

Moreover, for two spinor fields $\psi, \varphi$ we introduce a complex vector
field $\langle C \psi, \nabla
\varphi \rangle$ defined by the formula
\begin{displaymath}
\langle C \psi , \nabla \varphi \rangle \ := \ \langle C(X^k, X^l) \psi, 
\nabla_{X_l} \varphi \rangle \cdot X_k \ . 
\end{displaymath}

{\bf Lemma 1.2:} {\it For any $\psi \in \Gamma (S)$, we have the equation}
\begin{equation}
\begin{array}{c}
\Big\langle C(X^k, X^l) \nabla_{X_k} \psi , \nabla_{X_l} \psi \Big\rangle \ 
= \ \mbox{div} \Big\langle C \psi , \nabla \psi \Big\rangle - 
\frac{1}{2} \Big|C \psi \Big|^2 \\
\mbox{}\\
+ \ \frac{1}{4} \cdot \Big\langle \psi, \Big( ( \nabla_{X_k} \mbox{Ric})(X^l) 
\cdot X^k - X^k \cdot (\nabla_{X_k} \mbox{Ric})(X^l) \Big) \nabla_{X_l} 
\psi \Big\rangle . 
\end{array}
\end{equation}

{\bf Proof:} Let $x \in M^n$ be any point and let $(X_1, \ldots , X_n)$ be
any orthonormal frame in a neighbourhood of the point such that $(\nabla X_k)_x
=0$ holds for $k=1, \ldots , n$. Then we have at $x \in M^n$ that

\begin{displaymath}
\begin{array}{l}
\hspace{-1mm}\langle C(X^k , X^l) \nabla_{X_k} \psi , \nabla_{X_l} \psi 
\rangle \ = \ - \langle \nabla_{X_k} \psi , C(X^k , X^l) \nabla_{X_l} \psi 
\rangle \\[1em]
\hspace{-1mm} = - X_k (\langle \psi, C(X^k , X^l)\nabla_{X_l} \psi \rangle)
+ \langle \psi, (\nabla_{X_k} C)(X^k , X^l)\nabla_{X_l} \psi \rangle 
 + \langle \psi , C(X^k, X^l) \nabla_{X_k} \nabla_{X_l} \psi \rangle \\[1em]
\hspace{-1mm} = X_k (\langle C(X^k , X^l) \psi , \nabla_{X_l} \psi 
\rangle ) + \langle \psi , (\nabla_{X_k} C)(X^k , X^l) \nabla_{X_l} \psi 
\rangle + \frac{1}{2} \langle \psi , C(X^k, X^l)C(X_k , X_l) \psi 
\rangle \\[1em]
\hspace{-1mm} = \ \mbox{div} \langle C \psi , \nabla \psi  \rangle + \langle 
\psi 
, (\nabla_{X_k} C)(X^k , X^l) \nabla_{X_l} \psi \rangle - \frac{1}{2} 
|C \psi |^2 \ ,
\end{array}
\end{displaymath}

and we obtain the following formula for the left-hand side of 
the expression (10)
\begin{displaymath}
\textstyle
\langle C(X^k, X^l) \nabla_{X^k} \psi , \nabla_{X_l} \psi \rangle =  
\mbox{div} \langle
C \psi , \nabla \psi \rangle - \frac{1}{2} |C \psi |^2 
+ \langle \psi , (\nabla_{X_k} C)(X^k , X^l) \nabla_{X_l} \psi \rangle . 
\end{displaymath} 
On the other hand, from (9) we obtain
$(\nabla_Z C)(X,Y) =  \frac{1}{4} \cdot X^k \cdot (\nabla_Z K)(X,Y)(X_k)$ 
and
\begin{displaymath}
(\nabla_{X_k} C)(X^k , X^l)\ = \ \frac{1}{4} \cdot g \Big( (\nabla_{X_k} K)(X^k , X^l)(X^i), X^j \Big) X_i \cdot X_j \ .
\end{displaymath}
The Bianchi identity
\begin{displaymath}
(\nabla_X K)(Y,Z) + (\nabla_Y K)(Z,X) + (\nabla_Z K)(X,Y)\ = \ 0
\end{displaymath}

implies the relation
\begin{displaymath}
g((\nabla_{X_k} K)(X^k , X)(Y),Z)\ = \ g((\nabla_Z \mbox{Ric})(Y) 
- (\nabla_Y \mbox{Ric})(Z),X)\ .
\end{displaymath}

The latter two equations yield
\begin{displaymath}
\begin{array}{lll}
(\nabla_{X_k}C)(X^k, X^l)&=& \frac{1}{4} g\Big((\nabla_{X^j} \mbox{Ric})(X^i) 
- (\nabla_{X^i} \mbox{Ric})(X^j), X^l\Big) X_i \cdot X_j \\[1em]
&=& \frac{1}{4} \Big(g\Big(X^i, (\nabla_{X^j} \mbox{Ric})(X^l)\Big)-g\Big(X^j, 
(\nabla_{X_i} \mbox{Ric})(X^l)\Big)\Big) X_i \cdot X_j \\[1em]
&=& \frac{1}{4} \Big((\nabla_{X_j} \mbox{Ric})(X^l) \cdot X^j - 
X^i \cdot (\nabla_{X_i} \mbox{Ric})(X^l)\Big) . 
\end{array}
\end{displaymath}

Inserting this formula we obtain $(10)$. \hfill $\Box$\\

In the following we use the Schr\"odinger-Lichnerowicz formula
\begin{equation}
\nabla^* \nabla \ = \ D^2 - {\textstyle \frac{1}{4}}R \ .
\end{equation}

The local expression of the Bochner Laplacian $\nabla^* \nabla$ is
\begin{equation}
\nabla^* \nabla \ = \ - \nabla_{X_k , X^k} \ = \ - \nabla_{X_k} \nabla_{X^k} +
\Gamma_k \mbox{}^{kl} \nabla_{X_l} \ , 
\end{equation}

where the Christoffel symbols $\Gamma_{ij}\mbox{}^k$ are defined by
$\nabla_{X_i} X_j =  \Gamma_{ij} \mbox{}^k X_k$. In the proof of the following lemma we also use the well known general
formulas
\begin{equation}
X^k \cdot \nabla_{X,X_k} \psi \ = \ \nabla_X D \psi \ ,
\end{equation}
\begin{equation}{\textstyle
X^k \cdot \nabla_{X_k, X} \psi \ = \ D\nabla_X \psi - X^k \cdot
\nabla_{\nabla_{X_k}X} \psi \ = \ \nabla_X D \psi + \frac{1}{2} \mbox{Ric}(X) \cdot \psi . }
\end{equation}

Moreover, for $\psi, \varphi \in \Gamma (S)$, let $\psi \varphi$ and 
$\langle \psi, \nabla \varphi \rangle$ be the complex vector fields on $M^n$
locally given by
\begin{displaymath}
\psi \varphi \ := \ i \cdot \langle \psi , X^k \cdot \varphi \rangle \cdot X_k
 \quad , \quad 
\langle \psi , \nabla \varphi \rangle \ := \ \langle \psi , \nabla_{X^k} 
\varphi \rangle \cdot X_k \ . 
\end{displaymath}

The vector field satisfies the relation

\begin{equation}
i \cdot \mbox{div} (\psi \varphi)\ = \ \langle D \psi , \varphi \rangle - 
\langle \psi , D \varphi \rangle .
\end{equation}

{\bf Lemma 1.3:} {\it Let $\psi$ be any spinor field. Then there is the
identity}
\begin{equation}\mbox{} \quad 
\begin{array}{l}
\Big\langle C(X^k, X^l) \nabla_{X^k} \psi , \nabla_{X_l}  \psi \Big\rangle 
 =  - \mbox{Re} \left(\Big\langle \mbox{Ric} (X^k) \nabla_{X_k} D \psi ,  
\psi  \Big\rangle \right)
+ \Big\langle \nabla_{\mbox{Ric} (X^k)} \psi,  \nabla_{X^k} \psi \Big\rangle 
\\[1em]
+ \ \frac{1}{4}
\Big|\mbox{Ric}\Big|^2 \cdot \Big|\psi\Big|^2 - \frac{1}{2} \Big|C \psi 
\Big|^2 +
 \Big|\nabla D \psi \Big|^2 - \frac{R}{4} \cdot \Big|\nabla \psi \Big|^2 - 
\Big|(D^2 - \frac{R}{4}) \psi \Big|^2 \\[1em]
+  \mbox{div} \Big(i\left(\nabla_{X_k} D \psi + \frac{1}{2} \mbox{Ric} (X_k) 
\cdot \psi\right)\left(\nabla_{X_k} \psi\right)+ \Big\langle C \psi ,  \nabla 
\psi \Big\rangle - \Big\langle (D^2 - \frac{R}{4})\psi,  \nabla \psi 
\Big\rangle \Big)  . 
\end{array}
\end{equation}

{\bf Proof:} Let $x \in M^n$ be any point and let $(X_1 , \ldots , X_n)$ be 
any orthonormal frame in a neighbourhood of $x$ such that $(\nabla X_k)_x
=0$ for $k=1, \ldots , n$. We use the notations
\begin{displaymath}
R_{ijkl} \ := \ g (K(X_i, X_j)(X_k), X_l) , \quad
R_{ij} \ := \ g(\mbox{Ric}(X_i), X_j)= R_{ik}\mbox{}^k \mbox{}_j \ . 
\end{displaymath}

Then, we have $\Gamma_{ij}\mbox{}^k =0$ at the point $x$ and\\

$(*)$ \hfill $\displaystyle 
R_{ijk}\mbox{}^l =X_i (\Gamma_{jk}\mbox{}^l) - X_j (\Gamma_{ik}\mbox{}^l) , 
\quad R_i\mbox{}^j =X_i (\Gamma_k\mbox{}^{kj}) - X_k (\Gamma_i\mbox{}^{kj})\ 
. $ \hfill \mbox{}\\

\newcommand{\xkl}{(X^k, X^l)}

Using this we calculate
\begin{displaymath}
\begin{array}{l}
 \Big\langle C \xkl \nabla_{X_k} \psi , \nabla_{X_l} \psi \Big\rangle \ = \  
\Big\langle \nabla_{X^k} \nabla_{X^l} \nabla_{X_k} \psi - \nabla_{X^l} \nabla_{X^k} \nabla_{X_k} \psi , \nabla_{X_l} \psi \Big\rangle  \\[1em]
= \ \Big\langle \nabla_{X^k} \nabla_{X^l} \nabla_{X_k} \psi , 
\nabla_{X_l} \psi \Big\rangle + \Big\langle \nabla_{X^l} ( - \nabla_{X^k} 
\nabla_{X_k} \psi), \nabla_{X_l} \psi \Big\rangle \\[1em]
= \ \Big\langle \nabla_{X^k}(\nabla_{X^l, X_k} \psi + \Gamma_{lk}\mbox{}^p 
\nabla_{X_p} \psi), \nabla_{X^l} \psi \Big\rangle 
+ \Big\langle \nabla_{X_l} (- \nabla_{X^k, X_k} \psi - \Gamma_k\mbox{}^{kp} 
\nabla_{X_p} \psi), \nabla_{X^l} \psi \Big\rangle \\[1em]
\hspace{-4mm}\stackrel{(11)(12)}{=} \ \Big\langle \nabla_{X^k} (C(X_l, X_k)\psi + 
\nabla_{X_k, X_l} \psi ) , \nabla_{X^l} \psi \Big\rangle + X_k (\Gamma_{lk}
\mbox{}^p) \Big\langle \nabla_{X_p} \psi, \nabla_{X^l} \psi \Big\rangle \\[1em]
\hspace{6mm} + \ \Big\langle \nabla_{X_l} ((D^2 - \frac{R}{4})\psi), 
\nabla_{X^l} \psi \Big\rangle - X_l (\Gamma_k\mbox{}^{kp}) \Big\langle 
\nabla_{X_p} \psi, \nabla_{X^l} \psi \Big\rangle \\[1em]
= \ \Big\langle \nabla_{X_k} (C(X^l , X^k) \psi), \nabla_{X_l} \psi 
\Big\rangle + \Big\langle \nabla_{X^k} \nabla_{X_k, X_l} \psi , 
\nabla_{X^l} \psi 
\Big\rangle + (X_k (\Gamma_l\mbox{}^{kp} ) - X_l (\Gamma_k\mbox{}^{kp})) 
\cdot 
\end{array}
\end{displaymath}
\begin{displaymath}
\begin{array}{l}
\hspace{6mm} \Big\langle 
\nabla_{X_p} \psi, \nabla_{X^l} \psi \Big\rangle + X_l \Big(\Big\langle 
(D^2 - \frac{R}{4}) \psi , \nabla_{X^l} \psi 
\Big\rangle\Big)
- \Big\langle (D^2 - \frac{R}{4})\psi , \nabla_{X_l} \nabla_{X^l} \psi 
\Big\rangle\\[1em]
\stackrel{(*)}{=} \ X_k \Big(\Big\langle C(X^l, X^k) \psi , \nabla_{X_l} 
\psi \Big\rangle\Big) - \Big\langle C(X^l, X^k) \psi , \nabla_{X_k} \nabla_{X_l} \psi \Big\rangle \\[1em]
\hspace{6mm}+ \ \Big\langle \nabla_{X^k} (\nabla_{X_k} \nabla_{X_l} \psi - 
\Gamma_{kl}\mbox{}^p \nabla_{X_p} \psi), \nabla_{X^l} \psi \Big\rangle 
- R_l\mbox{}^p \Big\langle
\nabla_{X_p} \psi, \nabla_{X^l} \psi \Big\rangle \\[1em]
\hspace{6mm} + \ \mbox{div} \Big\langle (D^2 - \frac{R}{4} )\psi , \nabla 
\psi \Big\rangle + \Big|(D^2 - \frac{R}{4} )\psi \Big|^2  \\[1em]
= \ - \ \mbox{div} \Big\langle C \psi , \nabla \psi \Big\rangle + \frac{1}{2} 
|C \psi |^2 - \Big\langle (D^2 - \frac{R}{4}) \nabla_{X_k} \psi , \nabla_{X^k} 
\psi \Big\rangle - X^k (\Gamma_{kl}\mbox{}^p )\Big\langle \nabla_{X_p} \psi 
, \nabla_{X^l} \psi 
\Big\rangle \\[1em]
\hspace{6mm}-\ \Big\langle \nabla_{\mbox{Ric} (X^k)} \psi , \nabla_{X_k} \psi \Big\rangle + \mbox{div} \Big\langle (D^2 - \frac{R}{4} )\psi , \nabla \psi 
\Big\rangle + \Big|(D^2-\frac{R}{4} )\psi \Big|^2 . 
\end{array}
\end{displaymath}

Hence, it holds that
\begin{displaymath}
\begin{array}{ll}(2*) &
\begin{array}{l}
\Big\langle C(X^k , X^l) \nabla_{X_k} \psi , \nabla_{X_l} \psi
\Big\rangle \ = \ \Big|(D^2 - \frac{R}{4} )\psi \Big|^2 + \frac{1}{2} 
\Big|C \psi \Big|^2 - \Big\langle \nabla_{\mbox{Ric} (X^k)} \psi , 
\nabla_{X_k} \psi \Big\rangle \\[1em]
+ \ X_k (\Gamma^{kpl}) \cdot \Big\langle \nabla_{X_p} \psi , \nabla_{X_l} \psi \Big\rangle + \frac{R}{4} \cdot \Big|\nabla \psi \Big|^2 - \Big\langle D^2 \nabla_{X_l} \psi , 
\nabla_{X^l}\psi \Big\rangle \\[1em] 
+ \ \mbox{div} \Big(\Big\langle (D^2 - \frac{R}{4} )\psi , \nabla \psi \Big\rangle - \Big\langle C \psi , \nabla \psi \Big\rangle \Big) \ . 
\end{array}
\end{array}
\end{displaymath}

Further, we have
\begin{displaymath}
\begin{array}{l}
\Big\langle D^2 \nabla_{X_l} \psi , \nabla_{X^l} \psi \Big\rangle \ = \ 
\Big\langle D(D\nabla_{X_l} \psi ), \nabla_{X^l} \psi \Big\rangle \\[1em]
\stackrel{(14)}{=} \ \Big\langle D(\nabla_{X_l} D \psi + \frac{1}{2} \mbox{Ric}(X_l) \cdot \psi + \Gamma_{kl}\mbox{}^p X^k \cdot \nabla_{X_p} \psi), \nabla_{X^l} \psi \Big\rangle  \\[1em]
= \ \Big\langle D(\nabla_{X_l} D \psi + \frac{1}{2} \mbox{Ric} (X_l) \cdot 
\psi), \nabla_{X^l} \psi \Big\rangle + X_q (\Gamma_{kl}\mbox{}^p)
\Big\langle X^q \cdot X^k \cdot \nabla_{X_p} \psi , \nabla_{X^l} \psi 
\Big\rangle \\[1em]
\stackrel{(15)}{=} \ i \cdot \mbox{div} \Big((\nabla_{X_l} D \psi + \frac{1}{2} \mbox{Ric} (X_l) \cdot \psi)(\nabla_{X^l} \psi)\Big)+ \Big\langle 
\nabla_{X_l} D \psi + \frac{1}{2} \mbox{Ric}(X_l)
\cdot \psi, D \nabla_{X_l} \psi \Big\rangle \\[1em]
\hspace{6mm} + \left( \frac{1}{2} (X_q (\Gamma_{kl}\mbox{}^p) - X_k (\Gamma_{ql}\mbox{}^p )) + \frac{1}{2} (X_q (\Gamma_{kl}\mbox{}^p) + X_k (\Gamma_{ql}\mbox{}^k )) \right) \Big\langle X^q \cdot X^k \cdot \nabla_{X_p} \psi , \nabla_{X^l} \psi \Big\rangle  \\[1em]
\end{array}
\end{displaymath}
\begin{displaymath}
\begin{array}{l}
\stackrel{(15),(*)}{=} i \cdot \mbox{div} \Big(( \nabla_{X_l} D \psi + 
\frac{1}{2} \mbox{Ric} (X_l) \cdot \psi )(\nabla_{X^l} \psi )\Big)\\[1em]
\hspace{6mm}+ \ \Big\langle  \nabla_{X_l} D \psi + \frac{1}{2} \mbox{Ric} 
(X_l) \cdot \psi , \nabla_{X^l} D \psi + \frac{1}{2} \mbox{Ric} (X^l) \cdot 
\psi \Big\rangle \\[1em]
\hspace{6mm} +\  \frac{1}{2} R_{qkl}\mbox{}^p \langle X^q \cdot X^k \cdot \nabla_{X_p} \psi, \nabla_{X^l} \psi \rangle - X_k (\Gamma^{klp})\langle \nabla_{X_p} \psi, \nabla_{X_l} \psi \rangle  \\[1em]
= \ i \cdot \mbox{div} \Big((\nabla_{X_l} D \psi + \frac{1}{2} \mbox{Ric}
(X_l) \cdot \psi)(\nabla_{X^l} \psi )\Big) + |\nabla D \psi |^2 
- \mbox{Re} \Big(\langle \mbox{Ric} (X^l) \nabla_{X_l} D \psi , 
\psi \rangle \Big) \\[1em]
\hspace{6mm} +\ \frac{1}{4} |\mbox{Ric}|^2 |\psi|^2 - \frac{1}{2} R^{plqk} 
\langle X_q \cdot X_k \cdot \nabla_{X_p} \psi, \nabla_{X^l} \psi \rangle 
+ X_k (\Gamma^{kpl}) \langle \nabla_{X_p} \psi, \nabla_{X_l} \psi \rangle 
\\[1em]
\stackrel{(9)}{=} \ i \cdot \mbox{div} \Big((\nabla_{X_l} D \psi + \frac{1}{2} \mbox{Ric}(X_l) \cdot \psi)(\nabla_{X^l} \psi)\Big)+ |\nabla D \psi|^2 + 
\frac{1}{4} |\mbox{Ric}|^2 \cdot |\psi|^2\\[1em]
\hspace{6mm} - \  \mbox{Re} \Big(\langle \mbox{Ric} (X^l) \nabla_{X_l} 
D \psi , \psi \rangle \Big) + X_k (\Gamma^{kpl}) \langle \nabla_{X_p} \psi , 
\nabla_{X_l} \psi  \rangle - 2 \Big\langle C(X^p , X^l) \nabla_{X_p} \psi , 
\nabla_{X_l} \psi \Big\rangle \ . 
\end{array}
\end{displaymath}

Inserting the latter equation into $(2*)$ we obtain $(16)$. \hfill $\Box$\\

Comparing the equations $(10)$ and $(16)$ we obtain immediately\\

{\bf Lemma 1.4:} {\it For any spinor field $\psi$, we have the identity}
\begin{equation}
\begin{array}{l}
\mbox{Re} \Big(\Big\langle \mbox{Ric} (X^k) \nabla_{X_k} D \psi , 
\psi \Big\rangle \Big) + \frac{1}{4}
\Big\langle \psi , \Big((\nabla_{X_k} \mbox{Ric})(X^l)  X^k 
- X^k  ( \nabla_{X_k} \mbox{Ric})(X^l)\Big) \nabla_{X_l} \psi 
\Big\rangle \\[1em]
\mbox{} \quad 
= \ \Big|\nabla D \psi \Big|^2 - \Big|(D^2 - \frac{R}{4} ) \psi \Big|^2 - 
\frac{R}{4}  \Big|\nabla \psi \Big|^2 + \frac{1}{4} \Big|\mbox{Ric}
\Big|^2   \Big|\psi \Big|^2 
+ \Big\langle \nabla_{\mbox{Ric} (X^k)} \psi, \nabla_{X_k} \psi \Big\rangle  \\[1em]
\hspace{6mm} + \ \mbox{div}
\Big(i(\nabla_{X_k} D \psi + \frac{1}{2} \mbox{Ric} (X_k) \psi)(\nabla_{X^k}
\psi) - \Big\langle (D^2 - \frac{R}{4})\psi , \nabla \psi \Big\rangle \Big) \
. 
\end{array}
\end{equation}

The following purely algebraic condition on the covariant derivative of the
Ricci tensor implies that the second term in formula $(17)$ vanishes. The proof is an easy computation using the relations
in the Clifford algebra. A thorough geometric discussion of this condition 
will be provided in Section $2$.\\

{\bf Lemma 1.5:} {\it If the covariant derivative of the Ricci tensor satisfies}
\begin{displaymath}
(\nabla_X\mbox{Ric})(Y) \ = \ (\nabla_Y\mbox{Ric})(X) \ ,
\end{displaymath}
{\it then, for any spinor field $\psi$ and any vector field} $Y$, {\it the 
Clifford product}
\begin{displaymath}
\Big((\nabla_{X_k} \mbox{Ric})(Y)  \cdot X^k 
- X^k \cdot ( \nabla_{X_k} \mbox{Ric})(Y)\Big) \cdot \psi \ = \ 0
\end{displaymath}
{\it vanishes.}\\

We thus obtain the following Weitzenb\"ock formula for the length 
$|Q^t \psi|^2$, which is fundamental for all our further considerations.\\

{\bf Theorem 1.6:} {\it Let $M^n$ be a Riemannian spin manifold and suppose that}
\begin{displaymath}
(\nabla_X\mbox{Ric})(Y) \ = \ (\nabla_Y\mbox{Ric})(X) \ .
\end{displaymath}
{\it Then, for any spinor field $\psi$, there exists a vector field $X_{\psi} 
\in \Gamma (TM^n)$ such that}
\begin{equation}
\begin{array}{l}
\Big|Q^t \psi \Big|^2 \ = \ \Big|\nabla D \psi \Big|^2 - \frac{1}{n} 
\cdot \Big|D^2 \psi \Big|^2 + t^2 \cdot \Big|\mbox{Ric} - \frac{R}{n} \Big|^2 
\cdot \Big|\psi \Big|^2 
+ 2t \cdot \Big( \frac{R}{n} \cdot \mbox{Re} (\langle D^2 \psi , \psi \rangle ) \\[1em]
\hspace{14mm} + \ \Big|(D^2 - \frac{R}{4} )\psi \Big|^2 + \frac{R}{4} \cdot 
\Big|\nabla \psi \Big|^2  - \Big|\nabla D \psi \Big|^2 - \frac{1}{4} 
\cdot \Big|\mbox{Ric} \Big|^2 \Big|\psi \Big|^2 \\[1em]
\hspace{14mm} - \  \Big\langle \nabla_{\mbox{Ric} (X^k)} \psi , \nabla_{X_k} 
\psi \Big \rangle + \mbox{div} (X_{\psi} ) \Big) \ .
\end{array}
\end{equation}

{\bf Proof:} The formula follows from $(8)$ and $(17)$ if one defines
 $X_{\psi}$ locally by
\begin{displaymath}
\textstyle
X_{\psi} \ := \  \mbox{Re} \Big( \Big\langle (D^2 - \frac{R}{4}) \psi, \nabla \psi \Big\rangle - i \Big(\nabla_{X_k} D \psi + \frac{1}{2} \mbox{Ric} (X_k) 
\cdot \psi\Big)\Big(\nabla_{X^k} \psi \Big) \Big) .   \quad \quad \quad
\quad \Box
\end{displaymath}

\section{A mini-max principle for the estimate of the eigenvalues}

In this section we assume that $M^n$ is compact, connected and
that the Ricci tensor satisfies the condition
\begin{displaymath}
(\nabla_X\mbox{Ric})(Y) \ = \ (\nabla_Y\mbox{Ric})(X) \ .
\end{displaymath}

It is an easy consequence of the Bianchi idendity that the scalar curvature
of the manifold must be constant. Then the tensor
\begin{displaymath}
T(X) \ := \ \frac{1}{n-2} \Big( \frac{R}{2(n-1)}\cdot X - 
\mbox{Ric}(X) \Big)
\end{displaymath}
has the same properties as the Ricci tensor. In dimension $n=3$ the manifold is conformally flat. 
If $n \ge 4$, we obtain the identity 
\begin{displaymath}
(\nabla_{X_k}W)(X,Y,X^k) \ = \ (n-3) \cdot \big((\nabla_{X}T)(Y) - 
(\nabla_{Y}T)(X)\big)
\end{displaymath}
by computing the divergence of the Weyl tensor $W$ (see \cite{YB}). Therefore, 
the manifold satisfies the mentioned condition for the Ricci tensor if and 
only if it has constant scalar curvature and a harmonic Weyl tensor. Moreover, these two properties are equivalent to the condition that the curvature tensor
is harmonic (see Chapter $16$ in \cite{Besse}). The following examples are known:
\begin{enumerate}
\item Local products of Einstein manifolds;
\item conformally flat manifolds with constant scalar curvature;
\item warped products $S^1 \times_{f^2} N^{n-1}$ of an Einstein manifold with
positive scalar curvature $R = 4(n-1)/n$ by $S^1$ (see \cite{Der1}, \cite{Der2}, \cite{Der3}), where the function 
$F := f^{n/2}$ is a positive, periodic solution of the differential equation
\begin{displaymath}
F^{''} \ - \ F^{1-\frac{4}{n}} \ = \ - \, F \ ;
\end{displaymath}
\item warped products over Riemann surfaces.
\end{enumerate}
We denote by $\kappa_1(x) \le \kappa_2(x) \le \ldots \le \kappa_n(x)$ the eigenvalues of the Ricci tensor at the point $x \in M^n$ and by $\kappa_0$ the 
minimum of $\kappa_1$. If $D \psi = \lambda \psi$ is an 
eigenspinor, then $(18)$ yields the inequality
\begin{displaymath}
\int_{M^n} \Big( \Big( \frac{n-1}{n} \lambda^4 - \frac{R}{4} \lambda^2
+2 t (\frac{R}{n} \lambda^2 - \frac{1}{4} | \mbox{Ric} |^2 ) + t^2
|\mbox{Ric} - \frac{R}{n} |^2 \Big) | \psi |^2 
- 2t \Big\langle \nabla_{\mbox{Ric}(X^k)} \psi, \nabla_{X_k} \psi 
\Big\rangle \Big) \, \ge \, 0 \, . 
\end{displaymath}
In case of an Einstein manifold we get back the inequality $(1)$. In
general, the Schr\"odinger-Lichnerowicz formula and the estimation
\begin{displaymath}
\kappa_0 | \nabla \psi|^2 \ \le \ \langle \nabla_{\mbox{Ric} (X^k)} \psi , 
\nabla_{X_k} \psi \rangle
\end{displaymath}
imply the inequality
\begin{displaymath}
\kappa_0 \cdot \int_{M^n} (\lambda^2 -  \frac{R}{4}) |\psi|^2 
\ \le \ \int_{M^n} \langle \nabla_{\mbox{Ric} (X^k)} \psi , \nabla_{X_k}
\psi \rangle 
\end{displaymath}
and, finally, for any $t \ge 0$ we obtain the condition
\begin{displaymath}
\lambda^2\Big(\lambda^2 - \frac{nR}{4(n-1)}\Big) + 2t \frac{n}{n-1}
\Big( \frac{R}{n} - \kappa_0 \Big)\Big(\lambda^2 - \frac{R}{4}\Big)+
\frac{n}{n-1}\max_{x \in M^n}\Big[ \Big(t^2 - \frac{t}{2}\Big) \Big|\mbox{Ric} - \frac{R}{n}\Big|^2 \Big]  \ \ge \ 0 \ .
\end{displaymath}
This is a min-max principle and can be used in order to estimate the 
eigenvalues of the Dirac operator from below. Of course, only parameters 
between $0 \le t \le 1/2$ are interesting. A similar result involving only
the scalar curvature was proved in \cite{FKIM}. For $\lambda = 0$ we immediately obtain the following result\\

{\bf Theorem 2.1:} {\it Let $M^n$ be a compact Riemannian spin manifold with 
harmonic curvature tensor. If} $\kappa_0$ {\it and} $|\mbox{Ric}|^2_0$ {\it denote the minimum of the eigenvalues and the length of the Ricci tensor, respectively, and if}
\begin{displaymath}
|\mbox{Ric}|^2_0 \ > \ R \cdot \kappa_0
\end{displaymath}
{\it holds, then there are no harmonic spinors.}\\

If the scalar curvature is positive, we know that $\lambda^2 \ge nR/4(n-1)$ and
the mini-max principle yields a better estimate only in case that the left-hand
side is negative for $\lambda^2 = nR/4(n-1)$ and some $t > 0$. This
condition is equivalent to
\begin{displaymath}
 | \mbox{Ric} |^2_0 \ > \ \frac{R}{n-1} (R- \kappa_0) \ ,
\end{displaymath}
where $\kappa_0$ and $|\mbox{Ric} |_0$ are the minimum of the eigenvalues
and the length of the Ricci tensor, respectively\\

{\bf Example 1:} The warped product $S^1 \times_{f^2} N^{n-1}$ of an 
Einstein manifold $N^{n-1}$ with positive scalar curvature $R = 4(n-1)/n$  
by  $S^1$ never satisfies the latter condition. Solving the 
differential equation $F^{''} - F^{1-\frac{4}{n}} = - F$ for $\, n=5,\, F'(0) = 0$ and $F(0) = 0.1$ we obtain, for example, $\kappa_0 = - \, 8.5$ and $|\mbox{Ric}|_0^2 = 2$. One series of examples which we can apply our inequality 
to consist of products $ (S^1 \times_{f^2} N^{n-1})  \times  \ldots  \times  
 (S^1 \times_{f^2} N^{n-1})$ with a sufficiently large number of factors. A second series are products $\Sigma^k  \times  (S^1 \times_{f^2} N^{n-1})$ by an 
Einstein manifold $\Sigma^k$ with sufficiently large scalar curvature. We 
describe the case of a two-dimensional sphere $\Sigma^2$ and 
a $4$-dimensional Einstein spin manifold $N^4$ with scalar curvature 
$R_{N} = 16/5$ in greater delail. Consider the positive, periodic solution $F = f^{5/2}$ of the 
differential equation $F'' - F^{1/5} = - F$ with initial values $F(0) = 0.1$ 
and $F'(0) = 0$. The Ricci tensor of the manifold $S^1 \times_{f^2} N^{4}$ 
has two eigenvalues
\begin{displaymath}
\kappa_1 \ = \ \frac{24}{25} \Big( \frac{F'}{F}\Big)^2 + \frac{8}{5} \Big( 1
- F^{- \frac{4}{5}} \Big), \quad \kappa_2 \ = \ \frac{1}{4} \Big( \frac{16}{5}
- \kappa_1 \Big) \ .
\end{displaymath} 
The multiplicity of $\kappa_1$ is one, the multiplicity of $\kappa_2$ is four,
the scalar curvature of the warped product equals $16/5$. Denote by 
$R_{\Sigma}$ the scalar curvature of the sphere $\Sigma^2$ and consider the 
manifold $M^{7} := \Sigma^2 \, \times \, ( S^1 \times_{f^2} N^{4})$. Then we 
have
\begin{displaymath}
|\mbox{Ric}_M |_0^2 \ = \ \frac{R^2_{\Sigma}}{2} + |\mbox{Ric}_{(S^1 \times_{f^2} N^{4})}|_0^2 \ = \ 
\frac{R^2_{\Sigma}}{2} + 2, \quad R_M \ = \ R_{\Sigma} + \frac{16}{5}, \quad 
\kappa_0 \ = \ - \, 8.5 \ .
\end{displaymath}
For the optimal parameter $t= 0.212$ we obtain the estimate
\begin{displaymath}
\lambda^2 \ \ge \ - \, 1.74873 + 0.1105 \cdot R_{\Sigma} + 0.235194 \cdot
\sqrt{120.053 + 12.8828 \cdot R_{\Sigma} + R^2_{\Sigma}} \ \approx \ 
0.3457 \cdot R_{\Sigma} \ ,
\end{displaymath}
whereas the inequality $(1)$ yields the estimate $\lambda^2 \ge 
\frac{7}{24}R_{\Sigma} + \frac{14}{15} \approx 0.29 \cdot R_{\Sigma}$. $\hfill \Box$ \\

The discussion of the limiting case yields a spinor field $\psi$ in the kernel of one of the operators $Q^t$ with $t \ge 0$. Moreover, at every point we have
\begin{displaymath}
\kappa_1 | \nabla \psi|^2 \ = \ \langle \nabla_{\mbox{Ric} (X^k)} \psi , 
\nabla_{X_k} \psi \rangle \ ,
\end{displaymath}
i.e., the derivative of $\psi$ vanish in all directions $Y$ that are 
orthogonal to the $\kappa_1$-eigenspace of the Ricci tensor, 
$\nabla_Y\psi = 0$. The equation $Q^t \psi =0$ means that the eigenspinor 
$\psi$ satisfies the equation
\begin{displaymath}
\nabla_X \psi + \frac{\lambda}{n} X \cdot \psi + \frac{t}{\lambda}
(\mbox{Ric} - \frac{R}{n})(X) \cdot \psi \ =\ 0
\end{displaymath}
for each vector field $X$. In particular, the length of the spinor field is 
constant. If $t=0$, then $\psi$ is a Killing spinor. In case
$t > 0$, we consider the largest eigenvalue $\kappa_n$ at a minimum 
$x_0 \in M^n$ of $\kappa_1$ and insert an eigenvector. Then we obtain
\begin{displaymath}
\textstyle
\lambda^2 + t \cdot (n \cdot \kappa_n(x_0) - R) \ = \ 0  \ .
\end{displaymath}
But $n \cdot \kappa_n(x_0)-R$ is positive, a contradiction. Thus the limiting 
case in the inequality cannot occur except that $M^n$ is an Einstein manifold with a Killing spinor.\\

First we consider the case that the scalar curvature $R=0$ vanishes. Then 
$\kappa_0$ is negative and for any positive $t$ we have
\begin{displaymath}
\frac{n-1}{n} \lambda^4 - 2 t \cdot \kappa_0 \cdot\lambda^2 + \max_{x \in 
M^n}\Big[ \Big(t^2 - \frac{t}{2}\Big) |\mbox{Ric}|^2 \Big] \ > \ 0 \ .
\end{displaymath}
An elementary discussion yields the proof of the following theorem.\\

{\bf Theorem 2.2:} {\it Let $M^n$ be a compact, non-Ricci flat Riemannian 
spin manifold with harmonic curvature tensor and vanishing scalar curvature. If} $\kappa_0$ {\it and} $|\mbox{Ric}|^2_0$ {\it denote the minimum of the eigenvalues and the length of the Ricci tensor, respectively, then the eigenvalues of the Dirac operator are bounded by}
\begin{displaymath}
\lambda^2 > \frac{1}{4} \cdot \frac{| \mbox{Ric}|^2_0}{| \mbox{Ric}|_0 
\sqrt{\frac{n-1}{n}} + |\kappa_0|} \ . 
\end{displaymath} 

{\bf Remark:} The Schr\"odinger-Lichnerowicz formula implies the well 
known fact that a compact, non Ricci-flat Riemannian spin manifold 
with $R \equiv 0$ does not admit harmonic spinors. The estimate in 
Theorem 2.2 is a quantitative improvement of this fact for 
manifolds with harmonic curvature tensor and vanishing scalar curvature.\\

{\bf Example 2:} Let $\Gamma \subset \mbox{Conf}(S^n)$ be a geometrically finite
Kleinian group of compact type and denote by $\Lambda(\Gamma)$ its limit
set. Then $X^n(\Gamma) := \big(S^n - \Lambda(\Gamma)\big)/\Gamma$ is a closed
manifold equipped with a flat conformal structure. If the Hausdorff dimension 
of the limit set equals $(n-2)/2$, then there exists a Riemannian metric 
in the conformal class with vanishing scalar curvature (see \cite{SY}). 
S. Nayatani constructed this metric explicitly and studied its Ricci tensor 
(\cite{Nay}). \\

{\bf Example 3:} Let us continue Example 1. If $\Sigma^2$ is a compact surface
with scalar curvature $R_{\Sigma} = - \, \frac{16}{5}$, then 
$M^{7} := \Sigma^2 \, \times \, ( S^1 \times_{f^2} N^{4})$ has a harmonic 
curvature tensor and vanishing scalar curvature. Theorem 2.2 proves the 
estimate
\begin{displaymath}
\lambda^2 \ \ge \ 0.17833 \ .
\end{displaymath}

{\bf Example 4:}  If $\Sigma^2$ is a compact surface
with scalar curvature $R_{\Sigma} = - \, 4$, then 
$M^{7} := \Sigma^2 \, \times \, ( S^1 \times_{f^2} N^{4})$ has a harmonic 
curvature tensor and {\it negative}  scalar curvature. We apply the mini-max principle and obtain
\begin{displaymath}
\lambda^2 \ \ge \ 0.052 \ .
\end{displaymath}
In particular, $M^{7}$ has no harmonic spinors.
\section{An estimate of the eigenvalues}

Let us introduce the short cuts
\begin{displaymath}\textstyle
a \ := \ \frac{nR}{8(n-1)} , \quad b \ := \ \frac{n}{n-1} ( \frac{R}{n}- 
\kappa_0) , \quad c \ := \ |\mbox{Ric} - \frac{R}{n}|_0 \sqrt{\frac{n}{n-1}} ,\quad A \ := \ \textstyle \frac{c^2}{4} + 2 \frac{n-1}{n} ab 
\end{displaymath}
as well as the new parameter $s:= t/{\lambda^2}$. Then, by definition, $b$ and
$c$ are non-negative and the condition $|\mbox{Ric}|^2_0 \ > \ R \cdot \kappa_0$ of Theorem 2.1 is equivalent to $A > 0$. The mini-max principle yields 
immediately 
\begin{displaymath}
\lambda^2 \Big(\lambda^2 - 2 a + 2 b s (\lambda^2 - {\textstyle\frac{R}{4}}) +
(\lambda^2 s^2 - {\textstyle \frac{s}{2}}) c^2 \Big) \ \ge \ 0
\end{displaymath}
and then
\begin{displaymath}
\textstyle
\lambda^2 \ \ge \ \frac{2(a + A s)}{1 + 2 b s + c^2 s^2} \ =: \ f(s)
\end{displaymath}
for any $s \ge 0$. The function $f(s)$ attains its maximum at the point
\begin{displaymath}
\textstyle
s_0 \ := \ \frac{A - 2 a b}{a c^2 + c \sqrt{a^2 c^2 + A(A - 2 a b)}} \ .
\end{displaymath}
Hence, in case $R \le 0$ $(a \le 0)$, the parameter $s_0$ is automatically
positive. In case  $R > 0$ $(a > 0)$ the parameter $s_0 > 0$ is positive if and only if $A - 2 a b > 0$ or, equivalently, if
\begin{displaymath}
\textstyle
 | \mbox{Ric} |^2_0 \ > \ \frac{R}{n-1} (R- \kappa_0) \ ,
\end{displaymath}
holds. We summarize the main result.\\

{\bf Theorem 3.1:} {\it Let $M^n$ be a compact Riemannian spin manifold with
harmonic curvature tensor such that the condition}
\begin{equation}
\Big|\mbox{Ric} - \frac{R}{n} \Big|^2_0 \ > \ \Big( \frac{R}{n} - \kappa_0\Big) 
\max \left\{ 
\frac{R}{n-1} , - R \right\}
\end{equation}
{\it is satisfied. Then every eigenvalue $\lambda$ of the Dirac operator 
satisfies}
\begin{equation}
\lambda^2 \ > \ \frac{A^2}{bA - ac^2 + c \sqrt{a^2 c^2 + A(A - 2 a b)}} \ > 
\ 0 \ . 
\end{equation}

{\bf Proof:} In case $R \le 0$, the condition $(19)$ is equivalent to
$|\mbox{Ric}|^2_0 \ > \ R \kappa_0$ and the case of a positive scalar
curvature we already discussed. Then we obtain $\lambda^2 \ge f(s_0) $ and 
equality cannot occur for an eigenvalue $\lambda$ of $D$.  \hfill $\Box$\\

We remark that the compact, conformally flat $3$-manifolds with constant scalar curvature and constant length of the Ricci tensor are the $3$-dimensional 
space forms and the product of a $2$-dimensional space form $M^2$ by $S^1$ (see
\cite{CIS}). These manifolds do not satisfy the condition $(19)$. \\

In order to express the lower bound of the 
eigenvalue estimate in a convenient way we introduce the new variables 
$\alpha, \beta$ by the formulas
\begin{displaymath}
\alpha \ := \ a c^2 + (A - 2 a b)b, \quad 
\beta \ := \ (c^2 - b^2)(A - 2 a b)^2 \ .
\end{displaymath}
They are polynomials of degree three and six, depending on the eigenvalues of 
the Ricci tensor. The inequality $(20)$ can be reformulated in the following form.\\

{\bf Corollary 3.2:} {\it Let $M^n$ be a compact Riemannian spin manifold 
with harmonic curvature tensor and positive scalar curvature $R>0$. 
Suppose, moreover, that}
\begin{displaymath}
 \Big|\mbox{Ric} \Big|^2_0 \ > \ \frac{R}{n-1} \Big(R- \kappa_0\Big)
\end{displaymath}
{\it holds. Then each eigenvalue of the Dirac operator satisfies the estimate}
\begin{displaymath}
\lambda^2 \ > \ \frac{n R}{4 (n-1)} +\frac{(A - 2 ab)^2}
{\alpha + \sqrt{\alpha^2 + \beta}} \ > \  \frac{n R}{4 (n-1)} \ .
\end{displaymath}\\

We remark that the condition in Corollary 3.2 is satisfied in case that the
scalar curvature is positive and at least one eigenvalue of the Ricci tensor
is negative. Consequently, we obtain\\ 

{\bf Corollary 3.3:} {\it Let $M^n$ be a compact Riemannian spin manifold 
with parallel Ricci  tensor, positive scalar curvature and at least one negative eigenvalue of the Ricci tensor. Then the estimate of Corollary $3.2$ holds.}\\

{\bf Example 1:} Let us consider the product manifold $M^4= T^2 \times S^2$ 
equipped with the Riemannian metric induced by the metric of the flat torus
$T^2$ and the metric of the standard sphere $S^2 \subset {\Bbb R}^3$. Then
the Ricci tensor of $M^4$ is parallel and the set of eigenvalues of 
$\mbox{Ric}$ is given by  $(\kappa_1 , \ldots , \kappa_4)=(0,0,1,1)$. Hence, 
we have $\kappa_1=0, \, R=2 = |\mbox{Ric}|^2$. The condition $(19)$ is 
satisfied since
\begin{displaymath} \textstyle
|\mbox{Ric} - \frac{R}{4} |^2 \ = \ 1 \ > \ \frac{1}{3}
\ = \ \frac{1}{2} \max \left\{\frac{2}{3}, -2 \right\} \ . 
\end{displaymath}
Moreover, we find $a= \frac{1}{3} , \ b=\frac{2}{3} , \ c= \sqrt{\frac{4}{3}},
\ A= \frac{2}{3}, \ bA - ac^2 =0$. Inserting this into $(20)$ we obtain\\

$(*)$ \hfill $\lambda^2 > \frac{1}{2} \sqrt{2}$ \ . \hfill
\mbox{}\\

The Riemannian estimate $(1)$ yields the inequality $\lambda^2 \ge \frac{2}{3}$ and the K\"ahler estimate (see \cite{Ki}) gives the lower bound
$\lambda^2 \ge 1$. Since $\frac{2}{3} < \frac{1}{2}
\sqrt{2} <1$, the estimation of the first eigenvalue of the Dirac operator on the product considered as a K\"ahler manifold is the best one. We remark 
that $\lambda_1 = 1$ becomes an equality for the first 
eigenvalue (see \cite{F2}).\\

{\bf Example 2:} Let $N^2$ be any compact Riemannian surface of constant 
Gaussian curvature $-1$ and let $S^2 (r) \subset {\Bbb R}^3$ be the standard 
sphere of radius $r>0$. Then the Ricci tensor of the Riemannian product
$M^4 (r) := S^2 (r) \times N^2$ is parallel and $(\kappa_1 , \ldots , 
\kappa_4)=(-1 , -1, r^{-2} , r^{-2})$ is the corresponding
set of eigenvalues of Ric. Thus, here we have
\begin{displaymath} 
\textstyle
\kappa \ =\ -1, \quad R\ = \ - \, 2 \cdot (1- \frac{1}{r^2}), \quad 
|\mbox{Ric}|^2 \ =\ 2 \cdot (1+ \frac{1}{r^4}), \quad 
|\mbox{Ric} - \frac{R}{4} |^2 \ =\ (1+ \frac{1}{r^2})^2 \ ,
\end{displaymath}

and\\

\mbox{} \hspace{2cm} $
(\frac{R}{4} - \kappa) \cdot \max \left\{ \frac{R}{3} , \ - R\right\} \ 
= \ \max \left\{ \frac{1}{3} (\frac{1}{r^4} -1) , 1- \frac{1}{r^4} \right\} . $\\

This shows that the condition $(19)$ is satisfied and we find
\begin{displaymath} \textstyle
a\ = \ \frac{1}{3} ( -1 + \frac{1}{r^2} ) \quad , \quad b\ =\ \frac{2}{3}
( 1+ \frac{1}{r^2} ) \quad , \quad c\ = \ ( 1+ \frac{1}{r^2} )
\sqrt{\frac{4}{3}} , 
\end{displaymath}
\begin{displaymath} \textstyle
A\ = \ \frac{2}{3} \frac{1}{r^2} ( 1+ \frac{1}{r^2} ) \quad , \quad c^2-b^2
\ = \ \frac{8}{9} ( 1+ \frac{1}{r^2})^2 , \quad
bA - ac^2 \ = \ \frac{4}{9} ( 1+ \frac{1}{r^2} )^2 \ . 
\end{displaymath}
Inserting this into the inequality $(20)$ we obtain the estimation\\

$(*)$ \hfill $
\lambda^2 \ > \ \frac{1}{2} ( \sqrt{1+ \frac{2}{r^4}} - 1 ) \ > \ 0 \ . $
\hfill \mbox{}\\

Hence, we see that, on the product $M^4 (r)$, the Dirac operator has a trivial 
kernel even in case the scalar curvature is negative. In case the scalar curvature is positive,  we can compare the new estimation $(*)$ with the estimation 
$(1)$ and with the estimation for K\"ahler manifolds (see \cite{Ki}), 
respectively. For positive scalar curvature $(r <1)$, the lower bound $(*)$ is obviously better than the Riemannian estimate $(1)$,
\begin{displaymath}  \textstyle
\frac{2}{3} ( -1 + \frac{1}{r^2}) \ < \ \frac{1}{2} \Big( \sqrt{1+
\frac{2}{r^4}} -1 \Big) \ .
\end{displaymath}
If we compare the new lower bound $(*)$ and the lower bound $-1 + \frac{1}{r^2}$ in the K\"ahler case then, in the region  $1/\sqrt{2} <r<1 \ (0<R<2)$,  the inequality $(*)$ is the better one (see the figure):
\begin{displaymath} 
\textstyle
\frac{1}{2} \Big( \sqrt{1+\frac{2}{r^4}} -1 \Big) \ > \ -1 + \frac{1}{r^2} \ 
\quad \mbox{if} \quad 1/\sqrt{2} <r<1 \ .
\end{displaymath}
This example shows that in certain cases with positive scalar curvature the 
estimate given by Theorem 3.1 is even better than the bound in \cite{Ki} for K\"ahler manifolds. 
\begin{center}
\begin{displaymath}
\epsfig{figure=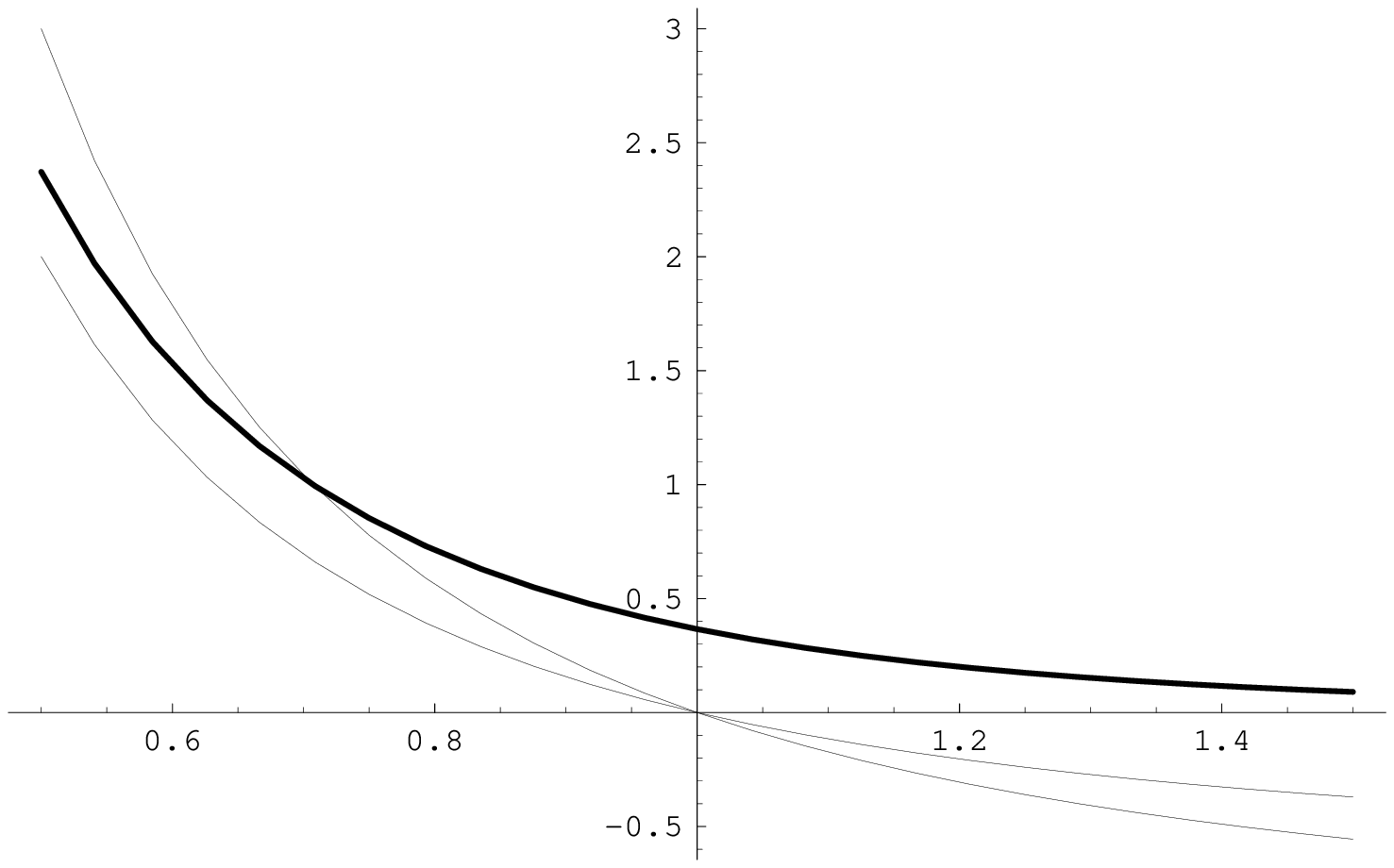, width=7cm}
\end{displaymath}
\end{center}
The preceding two examples are special cases of a more general 
situation. Consider compact Einstein manifolds with spin structures 
$M^{n_1}_1 , \ldots , M^{n_k}_k$ of dimensions $n_1,  \ldots , n_k $
(in the case of $n_i = 2$, we assume that $M^{n_i}$ is a surface of constant
Gaussian curvature). Then the Riemannian product $M^n := M^{n_1}_1 \times 
\ldots \times M^{n_k}_k$ is a compact Riemannian spin manifold with parallel 
Ricci tensor. Let $R_i$ be the scalar curvature of $M^{n_i}_i$. Then the 
scalar curvature $R$ as well as the length of the Ricci tensor of $M^n$ are 
given by
\begin{displaymath} 
\textstyle
R\ = \ \sum\limits^k_{i=1} R_i, \quad \quad 
|\mbox{Ric}|^2 \ = \ \sum\limits^k_{i=1} \frac{R_i^2}{n_i} \ . 
\end{displaymath}

Moreover, let us assume that the smallest eigenvalue of the Ricci tensor
is $\kappa_1 = \frac{R_1}{n_1}$. Then the conditions under which we can apply
our estimate are equivalent to
\begin{displaymath} 
\textstyle
\sum\limits^k_{i=1}
\frac{R_i^2}{n_i} \ > \ \frac{R_1}{n_1} ( \sum\limits^k_{i=1}
R_i ) \quad (R \ \le \ 0), \quad \mbox{and} \quad
 \sum\limits^k_{i=1}
\frac{R_i^2}{n_i} \ > \ \frac{1}{n-1} ( \sum\limits^k_{i=1} R_i )
( \sum\limits^k_{i=2} R_i ) \quad (R \ge 0) 
\end{displaymath}
 
respectively. Remark that, by the theorem of de Rham-Wu \cite{Wu},
 any compact, simply connected Riemannian manifold with parallel Ricci tensor 
splits into a Riemannian product of Einstein manifolds, i.e., the product situation is the general one for a parallel Ricci tensor. \\


\vspace{0.3cm}
Thomas Friedrich and Klaus-Dieter Kirchberg\\
Humboldt-Universit\"at zu Berlin\\
Institut f\"ur Reine Mathematik\\
Sitz: WBC Adlershof\\
D-10099 Berlin\\

\end{document}